\documentclass[11pt]{amsart}
\usepackage{amssymb}
\usepackage{amsmath}
\usepackage[active]{srcltx}
\usepackage{t1enc}
\usepackage[latin2]{inputenc}
\usepackage{verbatim}
\usepackage{amsmath,amsfonts,amssymb,amsthm}
\usepackage[mathcal]{eucal}
\usepackage{enumerate}
\usepackage[centertags]{amsmath}
\usepackage{graphicx}
\usepackage{color}
\graphicspath{ {./images/} }

\newtheorem{theorem}{Theorem}

\newtheorem{corollary}{Corollary}

\begin{document}
\author{N. Anakidze, N. Areshidze, L. Baramidze}

\title[Approximation by Nörlund means in Lebesgue spaces ]{Approximation by Nörlund means with respect to  Vilenkin system in Lebesgue spaces}
\address {N.Anakidze, The University of Georgia, School of Science and Technology, 77a Merab Kostava St, Tbilisi, 0171, Georgia.}
\email{nino.anakidze@mail.ru}
\address{N.Areshidze, Tbilisi State University, Faculty of Exact and Natural Sciences, Department of Mathematics, Chavchavadze str. 1, tbilisi 0128, georgia}
\email{nika.areshidze15@gmail.com}
\address{L.Baramidze, Tbilisi State University, Faculty of Exact and Natural Sciences, Department of Mathematics, Chavchavadze str. 1, tbilisi 0128, georgia}
\email{lashabara@gmail.com}

\thanks{}

%\date{}
\begin{abstract}
In this paper we improve and complement a result by M\'oricz and Siddiqi \cite{Mor}. In particular, we prove that their estimate of the Nörlund means with respect to the Vilenkin system holds also without their additional condition. Moreover, we prove a similar approximation result  in Lebesgue spaces for any $1\leq p<\infty$.
\end{abstract}
\maketitle 
\date{}



\section{Introduction}

Concerning some definitions and notations used in this introduction we refer to Section 2.   

It is well-known (see e.g. \cite{gol}, \cite{PTWbook} and \cite{Zy}) that,
 for any $1\leq p\leq \infty$ and $f\in L_p(G_m),$ there exists an absolute constant $C_p,$ depending only on $p$ such that
\begin{equation*}
\left\Vert \sigma_nf\right\Vert_p\leq C_p\left\|f\right\|_{p}.
\end{equation*}
Moreover, (for details see \cite{PTWbook})	if $1\leq p\leq \infty $, $M_N\leq n< M_{N+1}$, $f\in L^p(G_m)$ and $n\in \mathbb{N},$ then
\begin{eqnarray}\label{aaa}
\left\Vert \sigma_{n}f-f\right\Vert_{p}
\leq R^2\sum_{s=0}^{N} \frac{M_{s}}{M_{N}}\omega_p\left(1/M_s,f\right),
\end{eqnarray}
where $R:=\sup_{k\in\mathbb{N}}m_k$ and $\omega_p(\delta,f)$ is the modulus of continuity of $L^p, 1\le p \le\infty$  functions defined by
\begin{eqnarray*}
	\omega_p(\delta,f)=\sup_{|t|<\delta}\Vert f(x+t)-f(x) \Vert, \quad \quad \delta>0.
\end{eqnarray*}

%For each $\alpha>0$, Lipschitz classes in $L^p$ are defined by 

It follows that if $f\in lip\left( \alpha ,p\right) ,$ i. e.,

\begin{eqnarray*}
	Lip(\alpha,p):=\{f\in L^p: \omega_p(\delta,f)=O(\delta^{\alpha}) \quad\textrm{as}\quad \delta \to 0 \},
\end{eqnarray*}
then
	
	\begin{equation*}
	\left\Vert\sigma_nf-f\right\Vert _{p}=\left\{
	\begin{array}{c}O\left(1/M_N\right), \text{ \ \ \ if \ \ \ }\alpha >1,\\
	O\left(N/M_N\right),\text{ \ \ \ if \ \ \ }\alpha=1,\\
	O\left( 1/M_n^{\alpha}\right) ,\text{ \ \ \ \ \ \ \ if }\alpha <1.\end{array}\right.
	\end{equation*}
Moreover, (for details see \cite{PTWbook}) if $1\leq p< \infty ,$ $f\in L^{p}(G)$ and
	\begin{equation*}
	\left\Vert \sigma _{M_n}f-f\right\Vert_{p}=o\left( 1/M_n\right), \ \text{as} \ n\rightarrow \infty,
	\end{equation*}
then $f$  is a constant function.

The weak-$(1,1)$ type inequality  for the maximal operators of Vilenkin-Fejer means $\sigma ^{*}$, defined by
$$\sigma ^{*}f=\sup_{n\in\mathbb{N}}\vert \sigma_n f\vert$$ 
can be found in  Schipp
\cite{Sc} for Walsh series and in Pál, Simon \cite{PS} and Weisz \cite{We1}  for bounded Vilenkin
series. Boundedness of the maximal operators of Vilenkin-F\'ejer means of the one- and two-dimensional cases can be found in Fridli \cite{Fridli}, G\'at \cite{gat}, Goginava \cite{Goginava1}, Nagy and Tephnadze \cite{NT1,NT2,NT3,NT4}, Simon \cite{Simon1,Simon2}, Tutberidze \cite{tut1}, Weisz \cite{We3}.

Convergence and summability of  N\"orlund means with respect to Vilenkin systems were studied by Areshidze and Tephnadze\ \cite{AT1}, Blahota, Persson and Tephnadze \cite{BPT1} (see also \cite{BNPT1,BPTW1,BT2,BTT1}),  Fridli, Manchanda and  Siddiqi \cite{FMS}, Goginava \cite{Goginava}, Nagy \cite{na,nag,n} (see also  \cite{BN} and \cite{BNT}), Memic \cite{Memic}.

 M\'oricz and Siddiqi \cite{Mor} investigated the approximation properties of
some special N\"orlund means of Walsh-Fourier series of $L^{p}$ functions in
norm. In particular, they proved that if $f\in L^p(G),$ $1\leq p\leq \infty,$ $n=2^j+k,$ $1\leq k\leq 2^j \ (n\in \mathbb{N}_+)$ and $(q_k,k\in \mathbb{N})$ is a sequence of non-negative numbers, such that
\begin{equation}\label{MScond}
\frac{n^{\gamma-1}}{Q_n^{\gamma}}\sum_{k=0}^{n-1}q^{\gamma}_k =O(1),\ \ \text{for some} \ \ 1<\gamma\leq 2,
\end{equation}
then, there exists an absolute constant $C_p,$ depending only on $p$ such that

\begin{equation} \label{MSes}
\Vert t_nf-f\Vert_p\leq \frac{C_p}{Q_n} \sum_{i=0}^{j-1}2^iq_{n-2^i}\omega_p\left(\frac{1}{2^i},f\right)+C_p\omega_p\left( \frac{1}{2^j},f\right),
\end{equation}
when the sequence $(q_k,k\in \mathbb{N})$ is non-decreasing, while 

$$
\Vert t_nf-f\Vert_p\leq \frac{C_p}{Q_n} \sum_{i=0}^{j-1}\left(Q_{n-2^i+1}-Q_{n-2^{i+1}+1} \right)\omega_p\left(\frac{1}{2^i},f\right) +C_p\omega_p\left( \frac{1}{2^j},f\right),
$$
when the sequence  $(q_k,k\in \mathbb{N})$ is non-increasing.

In this paper we improve and complement a result by M\'oricz and Siddiqi \cite{Mor}. In particular, we prove that their estimate of the Nörlund means with respect to the Vilenkin system holds also without their additional condition. Moreover, we prove a similar approximation result  in Lebesgue spaces for any $1\leq p<\infty$.

\section{Preliminaries}

Let  $\mathbb{N}_{+}$ denote the set of the positive integers, $\mathbb{N}:=
\mathbb{N}_{+}\cup \{0\}.$ Let 
$m=:(m_0,m_1,...)$
be a sequence of positive integers not less than 2. Denote by
\begin{equation*}
Z_{m_k}:=\{0,1,...,m_k-1\}
\end{equation*}
the additive group of integers modulo $m_k.$ Define the group $G_m$ as the complete direct product of the group $%
Z_{m_k}$ with the product of the discrete topologies of $Z_{m_k}\textrm{'s}$.

The direct product $\mu $ of the measures 
\begin{equation*}
\mu_k\left( \{j\}\right):=1/m_k\text{ \ }(j\in Z_{m_k})  
\end{equation*}
is the Haar measure on $G_m$ with $\mu \left( G_m\right) =1.$

If $\sup_{k\in\mathbb{N}}m_k<+\infty$, then we call $G_m$ a bounded Vilenkin group. If $\{m_k\}_{k\ge0}$ sequance is unbounded, then $G_m$ is said to be unbounded Vilenkin group. In this paper we consider only bounded Vilenkin groups.

The elements of $G_m$ are represented by the sequences 

\begin{equation*}
x:=(x_{0},x_{1},\dots,x_{k},\dots)\qquad \left( \text{ }x_{k}\in
Z_{m_k}\right).
\end{equation*}
It is easy to give a base for the neighborhood of $G_m$, namely

\begin{eqnarray*}
I_{0}\left( x\right):=G_m, \ \ \
I_{n}(x):=\{y\in G_m\mid y_{0}=x_{0},\dots,y_{n-1}=x_{n-1}\}\text{ }(x\in
G_m,\text{ }n\in \mathbb{N}).
\end{eqnarray*}
For the simplicity we also define by $I_n$ as $I_n:=I_n(0).$

Let us define a generalized number system based on $m$ in the following way:
\begin{equation*}
    M_0=:1,\quad M_{k+1}=:m_{k}M_k\quad (k\in\mathbb{N})
\end{equation*}
Then every $n\in \mathbb{N}$ can be uniquely expressed as

$$
n=\sum_{k=0}^{\infty }n_{j}M_j, \ \ \ \text{where } \ \ \ n_{j}\in Z_{m_j}  \ \ \ (j\in \mathbb{
N})$$ 
and only a finite number of $n_{j}`$s differ from zero. Let 

 $$|n|=:\max\{j\in\mathbb{N}, n_j\neq0\}.$$

In 1947  Vilenkin \cite{Vi,Vi1,Vi2} investigated a group $G_m$ and introduced the Vilenkin systems $\{{\psi}_j\}_{j=0}^{\infty}$ as
\begin{equation*}
\psi _{n}\left( x\right):=\prod_{k=0}^{\infty }r_{k}^{n_{k}}\left( x\right) 
\text{ \quad }\left( n\in \mathbb{N}\right).
\end{equation*}
where $r_k(x)$ are the generalized Rademacher functions defined by
\begin{equation*}
r_{k}( x):=\exp(2\pi ix_k/m_k), \text{ \ \ }\left(k\in \mathbb{N}\right).
\end{equation*}

These systems include as a special case the Walsh system when $m_k=2$ for any $k\in\mathbb{N}$.

The norms (or quasi-norms) of Lebesgue spaces $L_{p}(G_{m})$  are defined by 
\begin{equation*}
\left\Vert f\right\Vert _{p}^{p}:=\int_{G_{m}}\left\vert f\right\vert
^{p}d\mu.
\end{equation*}%

The Vilenkin system is orthonormal and complete in $L^{2}\left( G_m\right)
\,$ (for details see e.g. \cite{AVD} and \cite{sws}).

If $f\in L^{1}\left(G_m\right) $, we can define the Fourier
coefficients, the partial sums of the Fourier series, the Fejér means, the Dirichlet and Fejér kernels 
with respect to the Vilenkin system in
the usual manner:

\begin{eqnarray*}
	\widehat{f}\left( k\right) &:&=\int_{G_m}f\overline{\psi}_{k}d\mu ,\,%
	\text{\quad }\left(k\in \mathbb{N}\right) , \\
	S_{n}f &:&=\sum_{k=0}^{n-1}\widehat{f}\left( k\right) \psi_{k},\text{
		\quad }\left(n\in \mathbb{N}_{+},\text{ }S_{0}f:=0\right) , \\
	\sigma _{n}f &:&=\frac{1}{n}\sum_{k=1}^{n}S_{k}f,\text{ \quad \ \  }\left(n\in \mathbb{N}_{+}\right).\\
 D_{n}&:&=\sum_{k=0}^{n-1}\psi_{k},\text{ \quad \ \  }\left(n\in \mathbb{N}_{+}\right).\\
 K_n&:&=\frac{1}{n}\sum_{k=1}^{n}D_{k},\text{ \quad \ \  }\left(n\in \mathbb{N}_{+}\right).\
\end{eqnarray*}

Recall that  (for details see e.g. \cite{AVD} and \cite{PTWbook}),
\begin{equation}\label{dn2.3}
\quad \hspace*{0in}D_{M_{n}}\left( x\right) =\left\{
\begin{array}{l}
\text{ }M_n,\text{\thinspace \thinspace \thinspace  if\thinspace
	\thinspace }x\in I_{n}, \\
\text{ }0,\text{\thinspace \thinspace \thinspace \thinspace \thinspace \thinspace \thinspace \thinspace if
	\thinspace \thinspace }x\notin I_{n},
\end{array}
\right.  
\end{equation}

\begin{eqnarray}\label{dn2.4}
D_{M_n-j}\left( x \right)&=&D_{M_n}\left( x \right)-\overline{\psi}_{M_n-1}(-x)D_{j}(-x)\\
&=&D_{M_n}\left( x \right)-\psi_{M_n-1}(x)\overline{D_{j}}(x),\,\,\ 0\le j<M_n.\nonumber
\end{eqnarray}
\begin{equation} \label{fn5}
n\left\vert K_n\right\vert\leq
2\sum_{l=0}^{\vert n\vert } M_l \left\vert K_{M_l} \right\vert,
\end{equation}
 and

\begin{eqnarray} \label{fn40}
\int_{G_m} K_n (x)d\mu(x)=1,  \ \ \ \ \ 
\sup_{n\in\mathbb{N}}\int_{G_m}\left\vert K_n(x)\right\vert d\mu(x)\leq 2.
\end{eqnarray}
Moreover, if $n>t,$ $t,n\in \mathbb{N},$ then 
\begin{equation}\label{lemma2}
K_{M_n}\left(x\right)=\left\{ \begin{array}{ll}
\frac{M_t}{1-r_t(x)},& x\in I_t\backslash I_{t+1},\quad x-x_te_t\in I_n, \\
\frac{M_n+1}{2}, & x\in I_n, \\
0, & \text{otherwise.} \end{array} \right.
\end{equation}

The $n$-th N\"orlund  mean $t_n$ of the Vilenkin-Fourier series of a integrable function $f$ is defined by
\begin{equation} \label{1.2}
t_nf:=\frac{1}{Q_n}\overset{n}{\underset{k=1}{\sum }}q_{n-k}S_kf,
\end{equation}
where
\begin{equation*} 
Q_n:=\sum_{k=0}^{n-1}q_k.
\end{equation*}
Here $\{q_k:k\geq 0\}$ is a sequence of nonnegative numbers, where $q_0>0$ and
$$
\lim_{n\rightarrow \infty }Q_{n}=\infty .
$$
Then the summability method (\ref{1.2}) generated by $\{q_k:k\geq 0\}$ is regular if and only if (see \cite{moo})
\begin{equation*}
\underset{n\rightarrow \infty }{\lim }\frac{q_{n-1}}{Q_{n}}=0.  \label{1a11}
\end{equation*}
In this paper we investigate regular N\"orlund  means only.

It is well-known (for details see e.g. \cite{PTWbook}) that  every  N\"orlund summability method generated by non-increasing  sequence $(q_k,k\in \mathbb{N})$ is regular, but N\"orlund means generated by non-decreasing  sequence $(q_k,k\in \mathbb{N})$ is not always regular.

The representation
\begin{equation*}
t_nf\left(x\right)=\underset{G_m}{\int}f\left(t\right)F_n\left(x-t\right) d\mu\left(t\right)
\end{equation*}
play central roles in the sequel, where 
\begin{equation}\label{1.3T}
F_n=:\frac{1}{Q_n}\overset{n}{\underset{k=1}{\sum }}q_{n-k}D_k
\end{equation}
is called the kernels of the N\"orlund  means.

If we invoke Abel transformation we get the following identities: 
\begin{eqnarray}  \label{2b}
Q_n:=\overset{n-1}{\underset{j=0}{\sum}}q_j=\overset{n}{\underset{j=1}{%
		\sum }}q_{n-j}\cdot 1 =\overset{n-1}{\underset{j=1}{\sum}}%
\left(q_{n-j}-q_{n-j-1}\right) j+q_0n
\end{eqnarray}
and
\begin{equation}  \label{2bbb}
t_nf=\frac{1}{Q_n}\left(\overset{n-1}{\underset{j=1}{\sum}}\left(
q_{n-j}-q_{n-j-1}\right) j\sigma_{j}f+q_0n\sigma_nf\right).
\end{equation}

\section{Formulation of Main Results}
Based on estimate \eqref{aaa} we can prove our next main results:
\begin{theorem}\label{Corollary3nnconv} 
	Let $M_N\leq n< M_{N+1}$ and ${{t}_{n}}$ be a regular N\"orlund  mean generated by non-decreasing sequence $\{q_k:k\in \mathbb{N}\},$ in sign $q_k \uparrow.$ Then, for some $f\in L^p(G_m),$ where $1\leq p< \infty, $ 
		
\begin{equation*}
    \Vert t_nf-f\Vert_p\leq \frac{3R^3}{Q_n} \sum_{i=0}^{N-1}M_iq_{n-M_i}\omega_p\left(\frac{1}{M_i},f\right)+2R^3\omega_p\left( \frac{1}{M_N},f\right).
\end{equation*}
\end{theorem}

\begin{theorem}\label{Corollary3nnconv0} 
	Let   ${{t}_{n}}$ be N\"orlund  mean generated by non-increasing sequence $\{q_k:k\in \mathbb{N}\}$, in sign $q_k \downarrow$.
	Then, for some $f\in L^p(G_m),$ where $1\leq p< \infty, $ 	
	
	\begin{eqnarray*}
		\Vert t_{M_n}f-f\Vert_p\leq 3R^2\sum_{s=0}^{n}\frac{M_s}{M_n}\omega _p\left(1/M_s,f\right)+C\overset{n-1}{\underset{s=0}{\sum}}\frac{(n-s)M_s}{M_{n}}\frac{q_{M_{s}}}{q_{M_{n}}}\omega _p\left(1/M_s,f\right).
	\end{eqnarray*}
\end{theorem}
\begin{theorem}\label{Corollary3nnconv1} 
Let $M_N\leq n< M_{N+1}$ and  ${{t}_{n}}$ be N\"orlund  mean generated by non-increasing sequence $\{q_k:k\in \mathbb{N}\}$, in sign $q_k \downarrow$, satisfying the condition

\begin{equation}\label{Cond}
\frac{1}{Q_n}=O\left(\frac{1}{n}\right),\text{\ \ as \ \ }n\rightarrow \infty.
\end{equation}
Then, for some $f\in L^p(G_m),$ where $1\leq p< \infty, $ 
	
\begin{eqnarray*}
\Vert t_nf-f\Vert_p \leq C\sum_{j=0}^{N}\frac{M_j}{M_N}\omega_p\left(1/M_j,f\right).
\end{eqnarray*}
\end{theorem}
As a consequence we obtain the following similar result proved in M\'oricz and Siddiqi \cite{Mor}:
\begin{corollary}
Let $\{ q_k : k \geq 0\}$ be a sequence of non-negative numbers such
that in case $q_k \uparrow$ condition 

\begin{equation} \label{fn00}
\frac{q_{n-1}}{Q_n}=O\left(\frac{1}{n}\right),\text{\ \ as \ \ }n\rightarrow \infty.
\end{equation}
 is satisfied, while in case $q_k \downarrow$ condition \eqref{Cond} is satisfied. If $f\in Lip(\alpha, p)$ for some $\alpha > 0$ and $1 \leq p< \infty,$ then
 
\begin{equation}\label{dn2.30}
\Vert t_nf-f\Vert_p=\left\{
\begin{array}{l}
\text{ }O(n^{-\alpha}),\text{\qquad \qquad  if\qquad  }0<\alpha<1, \\
\text{ }O(n^{-1}\log n),\text{\qquad  if
	\quad \ }\alpha=1,\\
\text{ }O(n^{-1}),\text{\qquad \qquad   if \quad \ \
	}\alpha>1,
\end{array}
\right.  
\end{equation}
\end{corollary}
As a consequence we obtain the following similar result proved in M\'oricz and Siddiqi \cite{Mor}:
\begin{corollary} a) Let  ${{t}_{n}}$ be  N\"orlund  means generated by non-decreasing sequence $\{q_k:k\in \mathbb{N}\}$ satisfying regularity condition \eqref{fn00}. If $f\in Lip(\alpha, p)$ for some $\alpha > 0$ and $1 \leq p< \infty,$ then ${t}_{n}f$ converge to $f$ in $L_p(G_m)$ norm.
	
	b)	Let  ${{t}_{n}}$ be N\"orlund  mean generated by non-increasing sequence $\{q_k:k\in \mathbb{N}\}$ satisfying condition \eqref{Cond}.
	If $f\in Lip(\alpha, p)$ for some $\alpha > 0$ and $1 \leq p<\infty,$  then ${t}_{n}f$ converge to $f$ in $L_p(G_m)$ norm.
\end{corollary}
\section{Proofs}

\begin{proof}[Proof of Theorem 1.]
	 Let $M_N\leq n<M_{N+1}.$ Since  $t_n$ be regular N\"orlund  means generated by the sequence of non-decreasing sequence $\{q_k:k\in \mathbb{N}\},$  by combining \eqref{2b} and \eqref{2bbb} we can conclude that
	\begin{eqnarray*}
		&&\Vert t_nf(x)-f(x)\Vert_p \\
		&\leq&\frac{1}{Q_n}\left(\overset{n-1}{\underset{j=1}{\sum}}\left(q_{n-j}-q_{n-j-1}\right)j\Vert\sigma_jf(x)-f(x)\Vert_p+q_0n\Vert\sigma_nf(x)-f(x)
		\Vert_p\right)\\
		%&\leq&\frac{1}{Q_n}\overset{n-1}{\underset{j=0}{\sum}}\left(q_{n-j}-q_{n-j-1}\right)j\alpha_j+\frac{q_0n\alpha_n}{Q_n}\\
		&:=&I+II.
	\end{eqnarray*}
Furthermore,
\begin{eqnarray*}
I&=&\frac{1}{Q_n}\overset{M_N-1}{\underset{j=1}{\sum}}\left(q_{n-j}-q_{n-j-1}\right)j\Vert\sigma_jf(x)-f(x)\Vert_p\\
&+&\frac{1}{Q_n}\overset{n-1}{\underset{j=M_N}{\sum}}\left(q_{n-j}-q_{n-j-1}\right)j\Vert\sigma_jf(x)-f(x)\Vert_p:=I_1+I_2.
\end{eqnarray*}
Now we estimate each terms separately. By applying \eqref{aaa} for $I_1$ we can conclude that
\begin{eqnarray}\label{I_1}
I_1&\leq&\frac{R^2}{Q_n}\overset{N-1}{\underset{k=0}{\sum}}\overset{M_{k+1}-1}{\underset{j=M_k}{\sum}}\left(q_{n-j}-q_{n-j-1}\right)j  
\sum_{s=0}^{k} \frac{M_{s}}{M_{k}}\omega_p\left(1/M_s,f\right)\\ \nonumber
&\leq&\frac{R^2}{Q_n}\overset{N-1}{\underset{k=0}{\sum}}M_{k+1}\overset{M_{k+1}-1}{\underset{j=M_k}{\sum}}\left(q_{n-j}-q_{n-j-1}\right) 
\sum_{s=0}^{k} \frac{M_{s}}{M_{k}}\omega_p\left(1/M_s,f\right)\\\nonumber
&\leq&\frac{R^3}{Q_n}\overset{N-1}{\underset{k=0}{\sum}}\left(q_{n-M_k}-q_{n-M_{k+1}}\right) \sum_{s=0}^{k} M_{s}\omega_p\left(1/M_s,f\right)\\\nonumber
&\leq&\frac{R^3}{Q_n}\overset{N-1}{\underset{s=0}{\sum}}M_{s}\omega_p\left(1/M_s,f\right)\sum_{k=s}^{N-1} \left(q_{n-M_k}-q_{n-M_{k+1}}\right)\\\nonumber
&\leq&\frac{R^3}{Q_n}\overset{N-1}{\underset{s=0}{\sum}}M_{s}q_{n-M_s}\omega_p\left(1/M_s,f\right).
\end{eqnarray}

It is evident that
\begin{eqnarray}\label{I_2}
	\quad \quad I_2&\leq&\frac{R^2}{Q_n}\overset{n-1}{\underset{j=M_N}{\sum}}\left(q_{n-j}-q_{n-j-1}\right)j  
	\sum_{s=0}^{N} \frac{M_{s}}{M_{N}}\omega_p\left(1/M_s,f\right)\\\nonumber
&\leq&\frac{R^2M_{N+1}}{Q_n}\overset{n-1}{\underset{j=M_N}{\sum}}\left(q_{n-j}-q_{n-j-1}\right)
\sum_{s=0}^{N} \frac{M_{s}}{M_{N}}\omega_p\left(1/M_s,f\right)\\\nonumber
&\leq&\frac{R^3q_{n-M_N}}{Q_n}
\sum_{s=0}^{N}M_{s}\omega_p\left(1/M_s,f\right)\\ \notag
&\leq&\frac{R^3}{Q_n}\overset{N}{\underset{s=0}{\sum}}M_{s}q_{n-M_s}\omega_p\left(1/M_s,f\right).
\end{eqnarray}

For $II$ we have that
\begin{eqnarray*}
	II&\leq&\frac{q_0R^2M_{N+1}}{Q_n}  
	\sum_{s=0}^{N} \frac{M_{s}}{M_{N}}\omega_p\left(1/M_s,f\right)\\
	&\leq& \frac{R^3}{Q_n}  
	\sum_{s=0}^{N-1}M_{s}q_{n-M_s}\omega_p\left(1/M_s,f\right)+R^3\omega_p\left(1/M_N,f\right).
\end{eqnarray*}
	The proof is complete.	
\end{proof}

\begin{proof}[Proof of Theorem 2.]
	By using \eqref{dn2.4} we find that
	\begin{eqnarray} \label{1.21}
	t_{M_n}f=D_{M_n}\ast f-\frac{1}{Q_{M_n}}\overset{M_n-1}{\underset{k=0}{\sum }}q_{k}\left(\left( \psi_{M_n-1}\overline{D_{k}}\right)\ast f\right).
	\end{eqnarray}
	By using Abel transformation we get
	\begin{eqnarray} 	\label{2cc}
	t_{M_n}f&=&D_{M_n}\ast f-\frac{1}{Q_{M_n}}\overset{M_n-2}{\underset{j=0}{\sum}}\left(q_j-q_{j+1}\right) j((\psi_{M_n-1}\overline{K_j})\ast f)\\ \notag
	&-& \frac{1}{Q_{M_n}}q_{M_n-1}(M_n-1)(\psi_{M_n-1}\overline{K}_{M_n-1}\ast f)\\ \notag
	&=&D_{M_n}\ast f-\frac{1}{Q_{M_n}}\overset{M_n-2}{\underset{j=0}{\sum}}\left(q_j-q_{j+1}\right) j((\psi_{M_n-1}\overline{K_j})\ast f)\\ \notag
	&-& \frac{1}{Q_{M_n}}q_{M_n-1}M_n(\psi_{M_n-1}\overline{K}_{M_n}\ast f)\\ \notag
	&+& \frac{q_{M_n-1}}{Q_{M_n}}(\psi_{M_n-1}\overline{D}_{M_n}\ast f)
	\end{eqnarray}
	and
	
	\begin{eqnarray} 	\label{2c}
	&&t_{M_n}f(x)-f(x)\\ \notag
	&=&\int_{G_m}(f(x+t)-f(x))D_{M_n}(t)dt\\ \notag
	&-& \frac{1}{Q_{M_n}}\overset{M_n-2}{\underset{j=0}{\sum}}\left(q_j-q_{j+1}\right) j\int_{G_m}\left(f(x+t)-f(x)\right)\psi_{M_n-1}(t)\overline{K}_j(t)dt\\  \notag
	&-& \frac{1}{Q_{M_n}}q_{M_{n}-1}M_n\int_{G_m}(f(x+t)-f(x))\psi_{M_n-1}(t)\overline{K}_{M_n}(t)dt\\ \notag
	&+& \frac{q_{M_{n}-1}}{Q_{M_n}}\int_{G_m}(f(x+t)-f(x))\psi_{M_n-1}(t)\overline{D}_{M_n}(t)dt\\ \notag
	&=&I+II+III+IV.
	\end{eqnarray}
	By combining generalized Minkowski's inequality and   \eqref{dn2.3}  we find that
	$$\Vert I\Vert _p\leq \int_{I_n}\Vert f(x+t)-f(x))\Vert_p D_{M_n}(t)dt\leq \omega_p\left(1/M_n,f\right). $$
and
	$$\Vert IV\Vert _p\leq \int_{I_n}\Vert f(x+t)-f(x))\Vert_p {D}_{M_n}(t)dt\leq \omega_p\left(1/M_n,f\right) .$$
 Since
$ M_nq_{M_n-1}\leq Q_{M_n} ,$ for any $n\in \mathbb{N},$ if combine \eqref{lemma2} and generalized Minkowski's inequality we get	
	\begin{eqnarray*}\label{fejaprox2}
		\Vert III\Vert_p
		&\leq&\int_{G_{m}}\left\Vert f\left( x+t\right) -f\left( x\right) \right\Vert_p
		\left|\overline{K}_{M_n}\left(t\right)\right| d\mu (t)\\ \notag
		&=&\int_{I_{n}}\left\Vert f\left( x+t\right) -f\left( x\right) \right\Vert_p \left|\overline{K}_{M_n}\left(t\right)\right|d\mu(t)\\
		&+&\sum_{s=0}^{n-1}\sum_{n_s=1}^{m_s-1}\int_{I_{n}\left(n_s e_{s}\right)}\left\Vert f\left( x+t\right) -f\left( x\right) \right\Vert_p \left|\overline{K}_{M_n}\left(t\right)\right| d\mu (t)\\ \notag
		&\leq&\int_{I_{n}}\left\Vert f\left( x+t\right)-f\left( x\right) \right\Vert_p\frac{M_{n}+1}{2}d\mu(t)\\ \notag
		&+&\sum_{s=0}^{n-1}M_{s+1}\sum_{n_s=1}^{m_s-1}\int_{I_n\left(n_se_s\right) }\left\Vert f\left(x+t\right)-f\left(x\right) \right\Vert_pd\mu(t)\\ \notag
		\\ \notag
		&\leq& \omega _{p}\left( 1/M_{n},f\right) \int_{I_{n}}\frac{M_{n}+1}{2}d\mu(t)\\\notag
		&+&\sum_{s=0}^{n-1}M_{s+1}\sum_{n_s=1}^{m_s-1}\int_{I_n\left(n_se_s\right)}\omega _{p}\left( 1/M_s,f\right)d\mu(t)\\ \notag
		&\leq&\omega _{p}\left( 1/M_{n},f\right)+R^2\sum_{s=0}^{n-1}\frac{M_s}{M_n}\omega _p\left(1/M_s,f\right)
		\leq R^2\sum_{s=0}^{n}\frac{M_s}{M_n}\omega _p\left(1/M_s,f\right).
	\end{eqnarray*}
 	This estimate also follows that	
 	
\begin{equation}
    \label{fejaprox2sssaaa}
         M_n\int_{G_m}\left\Vert f\left( x+t\right) -f\left( x\right) \right\Vert_p |\overline{K}_{M_{n}}(t)|\, d\mu (t)\leq R^2\sum_{s=0}^{n}M_s\omega _p\left(1/M_s,f\right).
\end{equation}

Let $M_k\leq j< M_{k+1}$ By applying  \eqref{fn5} and \eqref{fejaprox2sssaaa} we find that
	
\begin{eqnarray} \label{fejaprox22a}
        && j\int_{G_m}\left\Vert f\left( x+t\right) -f\left(x\right)\right\Vert_p
        |\overline{K}_{j}(t)|d\mu (t)\\ \notag
        &\leq& C\sum_{l=0}^{k}\sum_{s=0}^{l}M_s\omega _p\left(1/M_s,f\right).
\end{eqnarray}

	Hence, by combining \eqref{fn5} and \eqref{fejaprox22a} we find that
		
	\begin{eqnarray*}
		\Vert II\Vert_p&\leq& \frac{1}{Q_{M_n}}\overset{M_n-1}{\underset{j=0}{\sum}}\left(q_j-q_{j+1}\right) j\int_{G_m}\Vert f(x+t)-f(x)\Vert_p \vert\overline{K}_j(t)\vert d\mu (t)\\  \notag
		&\leq&\frac{1}{Q_{M_n}}\overset{n-1}{\underset{k=0}{\sum}}\overset{M_{k+1}-1}{\underset{j=M_k}{\sum}}\left(q_j-q_{j+1}\right) j\int_{G_m}\Vert f(x+t)-f(x)\Vert_p \vert \overline{K}_j(t)\vert d\mu (t)\\  \notag
		&\leq & \frac{C}{Q_{M_n}}\overset{n-1}{\underset{k=0}{\sum}}\overset{M_{k+1}-1}{\underset{j=M_k}{\sum}}\left(q_j-q_{j+1}\right)\sum_{l=0}^{k}\sum_{s=0}^{l}M_s\omega _p\left(1/M_s,f\right)
		\\  \notag
		&\leq & \frac{C}{Q_{M_n}}\overset{n-1}{\underset{k=0}{\sum}}\left(q_{M_k}-q_{M_{k+1}}\right) \sum_{l=0}^{k}\sum_{s=0}^{l}M_s\omega _p\left(1/M_s,f\right)
		\\  \notag
		&\leq & \frac{C}{Q_{M_n}}\overset{n-1}{\underset{l=0}{\sum}}\sum_{k=l}^{n-1}\left(q_{M_k}-q_{M_{k+1}}\right)\sum_{s=0}^{l}M_s\omega _p\left(1/M_s,f\right)
		\\  \notag
		&\leq & \frac{C}{Q_{M_n}}\overset{n-1}{\underset{l=0}{\sum}}q_{M_l}\sum_{s=0}^{l}M_s\omega _p\left(1/M_s,f\right) \\
		&\leq&  \frac{C}{Q_{M_n}}\overset{n-1}{\underset{s=0}{\sum}}M_s\omega _p\left(1/M_s,f\right)\sum_{l=s}^{n-1}q_{M_l}
		\\ \notag
		&\leq & \frac{C}{Q_{M_n}}\overset{n-1}{\underset{s=0}{\sum}}M_s\omega _p\left(1/M_s,f\right)q_{M_{s}}(n-s)\\
		&\leq& C\overset{n-1}{\underset{s=0}{\sum}}\frac{(n-s)M_s}{M_n}\frac{q_{M_{s}}}{q_{M_{n}}}\omega _p\left(1/M_s,f\right).
	\end{eqnarray*}

	The proof is complete.	
\end{proof}

\begin{proof} [Proof of theorem 3.]
	Let $M_N\leq n<M_{N+1}.$ Since  $t_n$ be regular N\"orlund  means, generated by sequence of non-increasing numbers $\{q_k:k\in \mathbb{N}\}$  by combining \eqref{2b} and \eqref{2bbb}, we can conclude that	
	\begin{eqnarray*}
		&&\Vert t_nf(x)-f(x)\Vert_p \\
		&\leq&\frac{1}{Q_n}\left(\overset{n-1}{\underset{j=1}{\sum}}\left(q_{n-j-1}-q_{n-j}\right)j\Vert\sigma_jf(x)-f(x)\Vert_p+q_0n\Vert\sigma_nf(x)-f(x)
		\Vert_p\right)\\
		&:=&I+II.
	\end{eqnarray*}	
	Furthermore,
	\begin{eqnarray*}
		I&=&\frac{1}{Q_n}\overset{M_{N}-1}{\underset{j=1}{\sum}}\left(q_{n-j-1}-q_{n-j}\right)j\Vert\sigma_jf(x)-f(x)\Vert_p\\
		&+& \frac{1}{Q_n}\overset{n-1}{\underset{j=M_{N}}{\sum}}\left(q_{n-j-1}-q_{n-j}\right)j\Vert\sigma_jf(x)-f(x)\Vert_p\\
		&=&I_1+I_2.
\end{eqnarray*}
Analogously to \eqref{I_1} we get that
\begin{eqnarray*}
	I_1&\leq&\frac{R^3}{Q_n}\overset{N-1}{\underset{k=0}{\sum}}\left(q_{n-M_{k+1}}-q_{n-M_{k}}\right)
	\sum_{s=0}^{k} M_{s}\omega_p\left(1/M_s,f\right)\\
	&\leq&\frac{R^3}{Q_n}\sum_{s=0}^{N-1}M_{s}	\omega_p\left(1/M_s,f\right)\overset{N-1}{\underset{k=s}{\sum}}\left(q_{n-M_{k+1}}-q_{n-M_{k}}\right)\\
	&=&\frac{R^3}{Q_n}\sum_{s=0}^{N-1}M_{s}	\omega_p\left(1/M_s,f\right)(q_{n-M_{N}}-q_{n-M_s})\\
	&\leq&\frac{R^3q_{n-M_{N}}}{Q_n}\sum_{s=0}^{N-1}M_{s}	\omega_p\left(1/M_s,f\right)\leq\frac{R^3q_0}{Q_n}\sum_{s=0}^{N-1}M_{s}	\omega_p\left(1/M_s,f\right).
\end{eqnarray*}

Analogously to \eqref{I_2} we find that

\begin{eqnarray*}
I_2&\leq&\frac{R^2}{Q_n}\overset{n-1}{\underset{j=1}{\sum}}\left(q_{n-j-1}-q_{n-j}\right)j  \sum_{s=0}^{N} \frac{M_{s}}{M_{N}}\omega_p\left(1/M_s,f\right)\\
&=&\frac{R^2}{Q_n}\left(nq_0-Q_{n}\right)\sum_{s=0}^{N} \frac{M_{s}}{M_{N}}\omega_p\left(1/M_s,f\right)\\
&\leq &\frac{M_{N+1}R^2q_0}{Q_nM_N}\sum_{s=0}^{N} M_{s}\omega_p\left(1/M_s,f\right)\\
&\leq&\frac{R^3q_0}{Q_n}\sum_{s=0}^{N} M_{s}\omega_p\left(1/M_s,f\right).
\end{eqnarray*}
	
For $II$ we have that
	\begin{eqnarray*}
II&\leq&\frac{q_0R^2M_{N+1}}{Q_n}\sum_{s=0}^{N} \frac{M_{s}}{M_{N}}\omega_p\left(1/M_s,f\right)\\
&\leq& \frac{R^3q_0}{Q_n} \sum_{s=0}^{N}M_{s}\omega_p\left(1/M_s,f\right).
\end{eqnarray*}
	
Using \eqref{Cond} we obtain estimate above so the proof is complete.
\end{proof}

\end{document}